\documentclass[preprint,12pt]{elsarticle}

\usepackage{amssymb}
\usepackage{amsmath}
\usepackage{amsthm}
\usepackage[colorlinks]{hyperref}
\AtBeginDocument{%
  \hypersetup{%
    linkcolor=blue,%
    citecolor=green,%
  }%
}
\usepackage{cleveref}

\crefname{equation}{}{}

\newtheorem{theorem}{Theorem}[section]
\newtheorem{lemma}[theorem]{Lemma}

\newtheorem{corollary}[theorem]{Corollary}

\theoremstyle{definition}
\newtheorem{definition}[theorem]{Definition}

\newtheorem{notation}[theorem]{Notation}

\theoremstyle{remark}
\newtheorem{remark}[theorem]{Remark}

\numberwithin{equation}{section}

\journal{}

\begin{document}

\begin{frontmatter}



\title{Boundary Pointwise $C^{1,\alpha}$ and $C^{2,\alpha}$ Regularity for Fully Nonlinear Elliptic Equations \tnoteref{t1}}


\author[rvt]{Yuanyuan Lian}
\ead{lianyuanyuan@nwpu.edu.cn; lianyuanyuan.hthk@gmail.com}
\author[rvt]{Kai Zhang\corref{cor1}}
\ead{zhang\_kai@nwpu.edu.cn; zhangkaizfz@gmail.com}
\tnotetext[t1]{This research is supported by the National Natural Science Foundation of China (Grant No. 11701454), the Natural Science Basic Research Plan in Shaanxi Province of China (Program No. 2018JQ1039) and the Fundamental Research Funds for the Central Universities (Grant No. 31020170QD032).}

\cortext[cor1]{Corresponding author. ORCID: \href{https://orcid.org/0000-0002-1896-3206}{0000-0002-1896-3206}}

\address[rvt]{Department of Applied Mathematics, Northwestern Polytechnical University, Xi'an, Shaanxi, 710129, PR China}

\begin{abstract}
In this paper, we obtain the boundary pointwise $C^{1,\alpha}$ and $C^{2,\alpha}$ regularity for viscosity solutions of fully nonlinear elliptic equations. I.e., If $\partial \Omega$ is $C^{1,\alpha}$ (or $C^{2,\alpha}$) at $x_0\in \partial \Omega$, the solution is $C^{1,\alpha}$ (or $C^{2,\alpha}$) at $x_0$. Our results are new even for the Laplace equation. Moreover, our proofs are simple.
\end{abstract}

\begin{keyword}
Boundary regularity \sep Schauder estimate \sep Fully nonlinear elliptic equation \sep Viscosity solution

\MSC[2010] 35B65 \sep 35J25 \sep 35J60 \sep 35D40

\end{keyword}

\end{frontmatter}


\section{Introduction}\label{S1}
Since 1980s, the fully nonlinear elliptic equations have been studied extensively (see \cite{MR1351007} and \cite{MR1118699} and the references therein). For the investigation on boundary behavior, there are usually two ways. One is to study the boundary regularity for viscosity solutions. Flattening the curved boundary by a transformation is widely applied (e.g. \cite{MR3246039}). However, the lower order terms and variant coefficients arise inevitably. Moreover, only local estimates can be derived rather than pointwise estimates. Another way is to obtain a priori estimates first and then use the method of continuity to prove the existence of classical solutions. It often requires more smoothness on the boundary and the boundary value (e.g. \cite{MR701522}). In both cases, the proofs are usually complicated. We note that in \cite{MR2853528}, Ma and Wang also proved the boundary pointwise $C^{1,\alpha}$ regularity by a barrier argument and a complicated iteration procedure.

In this paper, we study the boundary regularity for viscosity solutions and prove the pointwise $C^{1,\alpha}$ and $C^{2,\alpha}$ estimates under the corresponding pointwise geometric conditions on $\partial \Omega$. Our results are new even for the Laplace equation and these geometric conditions are rather general. Moreover, the boundaries don't need to be flattened and the proofs are simple.

The perturbation and compactness techniques are adopted here. We use solutions with flat boundaries to approximate the solution and the error between them can be estimated by maximum principles. Then, we can obtain the necessary compactness for solutions (see \Cref{l-7}). This basic perturbation idea is inspired originally by \cite{MR1351007}. The application to boundary regularity is inspired by \cite{MR3780142}. Based on the compactness result, we can obtain the desired estimates at the boundary if the boundary is ``almost'' flat (see \Cref{l-4} and \Cref{l-6}). This compactness technique has been inspired by \cite{MR3246039} and \cite{Wang_Regularity}. Then in aid of the scaling, the estimates on curved boundaries can be derived easily and the perturbation is a matter of scaling in some sense. The treatment for the right hand term and the boundary value is similar.

In this paper, we use the standard notations and refer to \Cref{no1.1} for details. Before stating our main results, we introduce the following notions.
\begin{definition}\label{d-f}
Let $A\subset R^n$ be a bounded set and $f$ be a function defined on $A$. We say that $f$ is $C^{k,\alpha}$ ($k\geq 0$) at $x_0\in A$ or $f\in C^{k, \alpha}(x_0)$ if there exist a polynomial $P$ of degree $k$ and a constant $K$ such that
\begin{equation}\label{holder}
  |f(x)-P(x)|\leq K|x-x_0|^{k+\alpha},~~\forall~x\in A.
\end{equation}
There may exist multiple $P$ and $K$ (e.g. $A=B_1\cap R^{n-1}$). Then we take $P_0$ with
\begin{equation*}
  \|P_0\|=\min \left\{\|P\|\big | ~\exists K ~\mbox{such that}\cref{holder} ~\mbox{holds with}~P~\mbox{and}~K\right\},
\end{equation*}
where $\|P\|=\sum_{m=0}^{k}|D^m P(x_0)|$. Define
\begin{equation*}
  D^mf(x_0)=D^mP_0(x_0),
\end{equation*}
\begin{equation*}
[f]_{C^{k,\alpha}(x_0)}=\min \left\{K\big | \cref{holder} ~\mbox{holds with}~P_0~\mbox{and}~K\right\}
\end{equation*}
and
\begin{equation*}
\|f\|_{C^{k, \alpha}(x_0)}=\|P_0\|+[f]_{C^{k, \alpha}(x_0)}.
\end{equation*}
\end{definition}

Next, we give the definitions of the geometric conditions on the domain.
\begin{definition}\label{d-re} Let $\Omega$ be a bounded domain and $x_0\in \partial \Omega$. We say that $\partial \Omega$ is $C^{k,\alpha}$ ($k\geq 1$) at $x_0$ or $\partial \Omega\in C^{k,\alpha}(0)$ if there exist a coordinate system $\{x_1,...,x_n \}$, a polynomial $P(x')$ of degree $k$ and a constant $K$ such that $x_0=0$ in this coordinate system,
\begin{equation}\label{e-re}
B_{1} \cap \{(x',x_n)\big |x_n>P(x')+K|x'|^{k+\alpha}\} \subset B_{1}\cap \Omega
\end{equation}
and
\begin{equation}\label{e-re2}
B_{1} \cap \{(x',x_n)\big |x_n<P(x')-K|x'|^{k+\alpha}\} \subset B_{1}\cap \Omega^c.
\end{equation}
Then, define
\begin{equation*}
[\partial \Omega]_{C^{k, \alpha}(x_0)}=\inf \left\{K\big|
\mbox{\cref{e-re} and\cref{e-re2} hold for }K \right\}
\end{equation*}
and
\begin{equation*}
\|\partial \Omega\|_{C^{k, \alpha}(x_0)}=\|P\|+[\partial \Omega]_{C^{k, \alpha}(x_0)}.
\end{equation*}

In addition, we define
\begin{equation*}
\underset{B_r}{\mathrm{osc}}~\partial\Omega=\underset{x\in \partial \Omega\cap B_r}{\sup} x_n -\underset{x\in \partial \Omega\cap B_r}{\inf} x_n.
\end{equation*}
\end{definition}

\begin{remark}\label{r-0}
Throughout this paper, we always assume that $0\in \partial \Omega$ and study the boundary behavior at $0$. When we say that $\partial \Omega$ is $C^{k,\alpha}$ at $0$, it always indicates that\cref{e-re} and\cref{e-re2} hold. Furthermore, without loss of generality, we always assume that
\begin{equation*}
  P(0)=0~\mbox{and}~DP(0)=0.
\end{equation*}
\end{remark}

\begin{remark}\label{r-1}
In this definition, $\partial \Omega$ doesn't need to be the graph of a function near $x_0$. For example, let
\begin{equation*}
  \Omega=B(e_n, 1)\backslash \left\{(x',x_n)\big | x_n=|x'|^2/2, |x|\leq 1/2 \right\}.
\end{equation*}
Then $\partial \Omega$ is $C^{2,\alpha}$ at $0$ by the definition. We will prove that the solution is $C^{2,\alpha}$ at $0$. Hence, our results are new even for the Laplace equation.
\end{remark}

Since we consider the viscosity solutions, the standard notions and notations for viscosity solutions are used, such as $\bar{S}(\lambda ,\Lambda ,f)$, $\underline{S}(\lambda ,\Lambda ,f)$, $S(\lambda ,\Lambda ,f)$, $M^{+}(M,\lambda ,\Lambda )$, $M^{-}(M,\lambda ,\Lambda )$ etc. For the details, we refer to \cite{MR1351007}, \cite{MR1376656} and \cite{MR1118699}. Without loss of generality, we always assume that the fully nonlinear operator $F$ is uniformly elliptic with ellipticity constants $\lambda$ and $\Lambda$, and $F(0)=0$. We call a constant $C$ universal if it depends only on $n,\lambda$ and $\Lambda$.

We use the Einstein summation convention in this work, i.e., repeated indices are implicitly summed over.

Now, we state our main results. For the boundary pointwise $C^{1,\alpha}$ regularity, we have
\begin{theorem}\label{t-1}
Let $0<\alpha<\alpha_1$ where $\alpha_1$ is a universal constant (see \Cref{l-1}). Suppose that $\partial \Omega$ is $C^{1,\alpha}$ at $0$ and $u$ satisfies
\begin{equation*}
\left\{\begin{aligned}
&u\in S(\lambda,\Lambda,f)&& ~~\mbox{in}~~\Omega\cap B_1;\\
&u=g&& ~~\mbox{on}~~\partial \Omega\cap B_1,
\end{aligned}\right.
\end{equation*}
where $g\in C^{1,\alpha}(0)$ and $f\in L^{n}(\Omega\cap B_1)$ satisfies for some constant $K_f$
\begin{equation}\label{e-dini-f}
  \|f\|_{L^{n}(\Omega\cap B_r)}\leq K_f r^{\alpha}, ~~\forall~0<r<1.
\end{equation}

Then $u\in C^{1,\alpha}(0)$, i.e., there exists a linear polynomial $L$ such that
\begin{equation}\label{e.t1-1}
  |u(x)-L(x)|\leq C |x|^{1+\alpha}\left(\|u\|_{L^{\infty }(\Omega\cap B_1)}+K_f+\|g\|_{C^{1,\alpha}(0)}\right), ~\forall ~x\in \Omega\cap B_{r_1},
\end{equation}
and
\begin{equation}\label{e.t1-2}
|Du(0)|\leq C \left(\|u\|_{L^{\infty }(\Omega\cap B_1)}+K_f+\|g\|_{C^{1,\alpha}(0)}\right),
\end{equation}
where $C$ depends only on $n, \lambda, \Lambda$ and $\alpha$, and $r_1$ depends also on $[\partial \Omega]_{C^{1,\alpha}(0)}$.
\end{theorem}

\begin{remark}\label{r-8}
In \cite{MR2853528}, Ma and Wang only proved the boundary pointwise $C^{1,\tilde{\alpha}}$ regularity for some $\tilde{\alpha}$ with $0<\tilde{\alpha}\leq \min(\alpha,\alpha_1)$ since the Harnack inequality was used. For instance, for the Laplace equation, we can obtain the $C^{1,\tilde{\alpha}}$ regularity for any $0<\tilde{\alpha}<1$, which can not been inferred from \cite{MR2853528}.
\end{remark}

For the boundary pointwise $C^{2,\alpha}$ regularity, we have
\begin{theorem}\label{t-2}
Let $0<\alpha<\alpha_2$ where $\alpha_2$ is a universal constant (see \Cref{l-2}). Suppose that $\partial \Omega$ is $C^{2,\alpha}$ at $0$ and $u$ satisfies
\begin{equation*}
\left\{\begin{aligned}
&F(D^2u)=f&& ~~\mbox{in}~~\Omega\cap B_1;\\
&u=g&& ~~\mbox{on}~~\partial \Omega\cap B_1,
\end{aligned}\right.
\end{equation*}
where $g\in C^{2,\alpha}(0)$ and $f\in C^{\alpha}(0)$.

Then $u\in C^{2,\alpha}(0)$, i.e., there exists a quadratic polynomial $P$ such that
\begin{equation}\label{en2}
  |u(x)-P(x)|\leq C |x|^{2+\alpha}\left(\|u\|_{L^{\infty }(\Omega\cap B_1)}+\|f\|_{C^{\alpha}(0)}+\|g\|_{C^{1,\alpha}(0)}\right), ~~\forall ~x\in \Omega\cap B_{r_1},
\end{equation}
and
\begin{equation}\label{e.t2-2}
|Du(0)|+\|D^2u(0)\|\leq C\left(\|u\|_{L^{\infty }(\Omega\cap B_1)}+\|f\|_{C^{\alpha}(0)}+\|g\|_{C^{1,\alpha}(0)}\right),
\end{equation}
where $C$ depends only on $n, \lambda, \Lambda$ and $\alpha$, and $r_1$ depends also on $\|\partial \Omega\|_{C^{2,\alpha}(0)}$.
\end{theorem}
\begin{remark}\label{r-4}
Note that the convexity of $F$ is not needed here, which is different from the interior $C^{2,\alpha}$ regularity.
\end{remark}

In the next section, we prepare some preliminaries. In particular, we prove the compactness and the closedness for a family of viscosity solutions. We obtain the boundary $C^{1,\alpha}$ regularity in \Cref{S3} and the boundary $C^{2,\alpha}$ regularity in \Cref{S4}.

\begin{notation}\label{no1.1}
\begin{enumerate}~~\\
\item $\{e_i\}^{n}_{i=1}$: the standard basis of $R^n$, i.e., $e_i=(0,...0,\underset{i^{th}}{1},0,...0)$.
\item $x'=(x_1,x_2,...,x_{n-1})$ and $x=(x_1,...,x_n)=(x',x_n)$ .
\item $S^n$: the set of $n\times n$ symmetric matrices and $\|A\|=$ the spectral radius of $A$ for any $A\in S^n$.
\item $R^n_+=\{x\in R^n\big|x_n>0\}$.
\item $B_r(x_0)=\{x\in R^{n}\big| |x-x_0|<r\}$, $B_r=B_r(0)$, $B_r^+(x_0)=B_r(x_0)\cap R^n_+$ and $B_r^+=B^+_r(0)$.
\item $T_r(x_0)\ =\{(x',0)\in R^{n}\big| |x'-x_0'|<r\}$ \mbox{ and } $T_r=T_r(0)$.
\item $A^c$: the complement of $A$ and $\bar A $: the closure of $A$, $\forall A\subset R^n$.
\item $\Omega_r=\Omega\cap B_r$ and $(\partial\Omega)_r=\partial\Omega\cap B_r$.
\item $\varphi _i=D_i \varphi=\partial \varphi/\partial x _{i}$ and $D \varphi=(\varphi_1 ,...,\varphi_{n} )$. Similarly, $\varphi _{ij}=D_{ij}\varphi =\partial ^{2}\varphi/\partial x_{i}\partial x_{j}$ and $D^2 \varphi =\left(\varphi _{ij}\right)_{n\times n}$.
\end{enumerate}
\end{notation}

\section{Preliminaries}
In this section, we introduce two lemmas stating the $C^{1,\alpha}$ and $C^{2,\alpha}$ regularity on flat boundaries. We will use them to approximate the solutions on curved boundaries. In addition, we prove the compactness and closedness for a family of viscosity solutions.

The following lemma concerns the boundary $C^{1,\alpha}$ regularity. It was first proved by Krylov \cite{MR688919} and further simplified by Caffarelli (see \cite[Theorem 9.31]{MR1814364} and \cite[Theorem 4.28]{MR787227}).
\begin{lemma}\label{l-1}
Let $u$ satisfy
\begin{equation*}
\left\{\begin{aligned}
&u\in S(\lambda,\Lambda,0)&& ~~\mbox{in}~~B_1^+;\\
&u=0&& ~~\mbox{on}~~T_1.
\end{aligned}\right.
\end{equation*}

Then there exists a universal constant $0<\alpha_1<1$ such that $u\in C^{1,\alpha_1}(0)$ and for some constant $a$,
\begin{equation*}
  |u(x)-ax_n|\leq C_1 |x|^{1+\alpha_1}\|u\|_{L^{\infty }(B_1^+)}, ~~\forall ~x\in B_{1/2}^+
\end{equation*}
and
\begin{equation*}
  |a|\leq C_1\|u\|_{L^{\infty }(B_1^+)},
\end{equation*}
where $C_1$ is universal.
\end{lemma}

The next lemma concerns the boundary $C^{2,\alpha}$ regularity. We refer \cite[Lemma 4.1]{MR3246039} for a proof.
\begin{lemma}\label{l-2}
Let $u$ satisfy
\begin{equation*}
\left\{\begin{aligned}
&F(D^2u)=0&& ~~\mbox{in}~~B_1^+;\\
&u=0&& ~~\mbox{on}~~T_1.
\end{aligned}\right.
\end{equation*}

Then there exists a universal constant $0<\alpha_2<1$ such that $u\in C^{2,\alpha_2}(0)$ and for some constants $a$ and $b_{in} (1\leq i \leq n)$,
\begin{equation}\label{e.4-1}
  |u(x)-ax_n-b_{in}x_ix_n|\leq C_2 |x|^{2+\alpha_2}\|u\|_{L^{\infty }(B_1^+)}, ~~\forall ~x\in B_{1/2}^+,
\end{equation}
\begin{equation}\label{e.4-2}
  F(b_{in})=0
\end{equation}
and
\begin{equation*}
  |a|+|b_{in}|\leq C_2\|u\|_{L^{\infty }(B_1^+)},
\end{equation*}
where $C_2$ is universal.
\end{lemma}
\begin{remark}\label{r-6}
In\cref{e.4-1}, the Einstein summation convention is used (similarly hereinafter). In\cref{e.4-2}, $b_{in}$ denotes the matrix $a_{ij}$ whose elements are all $0$ except $a_{in}=b_{in}$ for $1\leq i\leq n$ (similarly hereinafter).
\end{remark}

The following lemma presents a uniform estimate for solutions, which is a kind of ``equicontinuity'' up to the boundary.
\begin{lemma}\label{l-3}
Let $0<\delta<1/4$. Suppose that $u$ satisfies
\begin{equation*}
\left\{\begin{aligned}
&u\in S(\lambda,\Lambda,f)&& ~~\mbox{in}~~\Omega_1;\\
&u=g&& ~~\mbox{on}~~(\partial \Omega)_1,
\end{aligned}\right.
\end{equation*}
with $\|u\|_{L^{\infty}(\Omega_1)}\leq 1$, $\|f\|_{L^{n}(\Omega_1)}\leq\delta$, $\|g\|_{L^{\infty}((\partial \Omega)_1)}\leq \delta$
and $\underset{B_1}{\mathrm{osc}}~\partial\Omega \leq \delta$.

Then
\begin{equation*}
  \|u\|_{L^{\infty}(\Omega_{\delta})}\leq C\delta,
\end{equation*}
where $C$ is universal.
\end{lemma}

\proof Let $\tilde{B}^{+}_{1}=B^{+}_{1}-\delta e_n $ and $\tilde{T}_1=T_1-\delta e_n$. Then $(\partial \Omega)_{1/4}\subset \tilde{B}^{+}_{1}$. Let $v$ solve
\begin{equation*}
\left\{\begin{aligned}
 &M^{+}(D^2v,\lambda,\Lambda)=0 &&\mbox{in}~~\tilde{B}^{+}_{1}; \\
 &v=0 &&\mbox{on}~~\tilde{T}_{1};\\
 &v=1 &&\mbox{on}~~\partial \tilde{B}^{+}_{1}\backslash \tilde{T}_{1}.
\end{aligned}
\right.
\end{equation*}
Let $w=u-v$ and then $w$ satisfies (note that $v\geq 0$)
\begin{equation*}
    \left\{
    \begin{aligned}
      &w\in \underline{S}(\lambda /n,\Lambda , f) &&\mbox{in}~~ \Omega \cap \tilde{B}^{+}_{1}; \\
      &w\leq g &&\mbox{on}~~\partial \Omega \cap \tilde{B}^{+}_{1};\\
      &w\leq 0 &&\mbox{on}~~\partial \tilde{B}^{+}_{1}\cap \bar{\Omega}.
    \end{aligned}
    \right.
\end{equation*}

By \Cref{l-1},
\begin{equation*}
  \|v\|_{L^{\infty}(\tilde{B}_{4\delta}^+)}\leq C\delta,
\end{equation*}
where $C$ is universal. For $w$, by the Alexandrov-Bakel'man-Pucci maximum principle, we have
\begin{equation*}
  \begin{aligned}
\sup_{\Omega \cap \tilde{B}^{+}_{1}} w& \leq\|g\|_{L^{\infty }(\partial \Omega \cap \tilde{B}^{+}_{1})}+C\|f\|_{L^n(\Omega \cap \tilde{B}^{+}_{1})}\leq C\delta,
  \end{aligned}
\end{equation*}
where $C$ is universal. Hence,
\begin{equation*}
\sup_{\Omega_{\delta}} u\leq \sup_{\Omega\cap \tilde{B}_{4\delta}^+} u \leq \|v\|_{L^{\infty}(\tilde{B}_{4\delta}^+)}+\sup_{\Omega \cap \tilde{B}^{+}_{1}} w\leq C\delta.
\end{equation*}

The proof for
\begin{equation*}
  \inf_{\Omega_{\delta}} u \geq -C\delta
\end{equation*}
is similar and we omit it here. Hence, the proof is completed.\qed~\\

\begin{remark}\label{r-2}
The proof shows the idea that approximating the general solution $u$ by a solution $v$ with a flat boundary. This idea is inspired by \cite{MR3780142}.
\end{remark}

Based on the above lemma, the following corollary follows easily:
\begin{corollary}\label{c-1}
For any $0<r<1$ and $\varepsilon>0$, there exists $\delta>0$ (depending only on $n,\lambda,\Lambda,r$ and $\varepsilon$) such that if $u$ satisfies
\begin{equation*}
\left\{\begin{aligned}
&u\in S(\lambda,\Lambda,f)&& ~~\mbox{in}~~\Omega_1;\\
&u=g&& ~~\mbox{on}~~(\partial \Omega)_1,
\end{aligned}\right.
\end{equation*}
with $\|u\|_{L^{\infty}(\Omega_1)}\leq 1$, $\|f\|_{L^{n}(\Omega_1)}\leq \delta$, $\|g\|_{L^{\infty}((\partial \Omega)_1)}\leq \delta$
and $\underset{B_1}{\mathrm{osc}}~\partial\Omega \leq \delta$, then
\begin{equation*}
  \|u\|_{L^{\infty}(\Omega \cap B(x_0,\delta))}\leq \varepsilon, ~~\forall~x_0\in \partial \Omega\cap B_r.
\end{equation*}
\end{corollary}

Next, we prove the equicontinuity of the solutions, which provides the necessary compactness.

\begin{lemma}\label{l-7}
For any $\Omega'\subset\subset \bar{\Omega}\cap B_1$ and $\varepsilon>0$, there exists $\delta>0$ (depending only on $n,\lambda,\Lambda,\Omega'$ and $\varepsilon$) such that if $u$ satisfies
\begin{equation*}
\left\{\begin{aligned}
&u\in S(\lambda,\Lambda,f)&& ~~\mbox{in}~~\Omega_1;\\
&u=g&& ~~\mbox{on}~~(\partial \Omega)_1,
\end{aligned}\right.
\end{equation*}
with $\|u\|_{L^{\infty}(\Omega_1)}\leq 1$, $\|f\|_{L^{n}(\Omega_1)}\leq \delta$, $\|g\|_{L^{\infty}((\partial \Omega)_1)}\leq \delta$
and $\underset{B_1}{\mathrm{osc}}~\partial\Omega \leq \delta$, then for any $x,y\in \Omega'$ with $|x-y|\leq \delta$, we have
\begin{equation*}
  |u(x)-u(y)|\leq \varepsilon.
\end{equation*}
\end{lemma}

\proof By \Cref{c-1}, for any $\varepsilon>0$, there exists $\delta_1>0$ depending only on $n,\lambda,\Lambda,\varepsilon$ and $\Omega'$ such that for any $x,y\in \Omega'$ with $\mathrm{dist}(x,\partial \Omega)\leq \delta_1$ and $|x-y|\leq \delta_1$, we have
\begin{equation}\label{e.l3-2.1}
|u(x)-u(y)|\leq |u(x)|+|u(y)|\leq \varepsilon.
\end{equation}

If $\mathrm{dist}(x,\partial \Omega)> \delta_1$, by the interior H\"{o}lder estimate,
\begin{equation}\label{e.l3-2.2}
|u(x)-u(y)|\leq C\frac{|x-y|^{\alpha}}{\delta_1^{\alpha}},
\end{equation}
where $C$ and $0<\alpha<1$ are universal. Take $\delta$ small enough such that
\begin{equation*}
  C\frac{\delta^{\alpha}}{\delta_1^{\alpha}}\leq \varepsilon.
\end{equation*}
Then by combining\cref{e.l3-2.1} and\cref{e.l3-2.2}, the conclusion follows.~\qed~\\

Now, we give a closedness result for viscosity solutions.

\begin{lemma}\label{l-5}
Let $u_k\in C(\bar{\Omega}_k\cap B_1)$ $(k\geq 1)$ satisfy
\begin{equation*}
\left\{\begin{aligned}
&F_k(D^2u_k)\geq (\leq ) f_k&& ~~\mbox{in}~~\Omega_k\cap B_1;\\
&u_k=g_k&& ~~\mbox{on}~~\partial \Omega_k\cap B_1.
\end{aligned}\right.
\end{equation*}

Suppose that $F_k\rightarrow F$ uniformly on compact subsets of $S^n$, $\|f_k\|_{L^{n}(\Omega_k\cap B_1)}\rightarrow 0$, $\|g_k\|_{L^{\infty}(\partial \Omega_k\cap B_1)}\rightarrow 0$ and $\underset{B_1}{\mathrm{osc}}~\partial\Omega_k \rightarrow 0$.

In addition, assume that for any $\Omega'\subset\subset\bar{\Omega}\cap B_1$, $u_k\rightarrow u$ uniformly on $\Omega'$. That is, for any $\varepsilon>0$, there exists $k_0$ such that for any $k\geq k_0$ and $x\in \Omega'\cap \bar{\Omega}_k$, we have
\begin{equation*}
  |u_k(x)-u(x)|\leq \varepsilon.
\end{equation*}

Then $u\in C(B_1^+\cup T_1)$ and
\begin{equation*}
\left\{\begin{aligned}
&F(D^2u)\geq (\leq) 0&& ~~\mbox{in}~~ B_1^+;\\
&u=0&& ~~\mbox{on}~~T_1,
\end{aligned}\right.
\end{equation*}
\end{lemma}

\proof We only prove the case for a subsolution. From \cite[Theorem 3.8]{MR1376656}, $F(D^2u)\geq 0$ in $B_1^+$ holds. For any $x_0\in T_1$ and $\varepsilon>0$, let $\delta>0$ be small to be specified later and $\tilde{x}\in B^+(x_0,\delta)\subset\subset B_1$. Since $u_k$ converges to $u$ uniformly, there exists $k_0$ such that for any $k\geq k_0$ and $x\in B^+(x_0,\delta)\cap \bar{\Omega}_k$, we have
\begin{equation*}
  |u_k(x)-u(x)|\leq \varepsilon/2.
\end{equation*}
Take $k$ large enough such that $\tilde{x}\in \Omega_k$ and $\|g_k\|_{L^{\infty}(\partial \Omega_k\cap B_1)}\leq \varepsilon/4$. Note that $u_k\in C(\bar{\Omega}_k\cap B_1)$. Then we can take $\delta$ small such that $|u_k(\tilde{x})|\leq \varepsilon/2$. Hence,
\begin{equation*}
  |u(\tilde{x})|=|u(\tilde{x})-u_k(\tilde{x})+u_k(\tilde{x})|\leq |u(\tilde{x})-u_k(\tilde{x})|+|u_k(\tilde{x})|\leq \varepsilon.
\end{equation*}
Therefore, $u$ is continuous up to $T_1$ and $u\equiv 0$ on $T_1$. ~\qed~\\

\section{Boundary $C^{1,\alpha}$ regularity}\label{S3}
In this section, we give the proof of the boundary $C^{1,\alpha}$ regularity. First, we prove that the solution in \Cref{t-1} can be approximated by a linear function provided that the prescribed data are small enough.

\begin{lemma}\label{l-4}
Let $\alpha_1$ and $C_1$ be as in \Cref{l-1}. For any $0<\alpha<\alpha_1$, there exists $\delta>0$ such that if $u$ satisfies
\begin{equation*}
\left\{\begin{aligned}
&u\in S(\lambda,\Lambda,f)&& ~~\mbox{in}~~\Omega_1;\\
&u=g&& ~~\mbox{on}~~(\partial \Omega)_1,
\end{aligned}\right.
\end{equation*}
with $\|u\|_{L^{\infty}(\Omega_1)}\leq 1$, $\|f\|_{L^{n}(\Omega_1)}\leq \delta$, $\|g\|_{L^{\infty}((\partial \Omega)_1)}\leq \delta$
and $\underset{B_1}{\mathrm{osc}}~\partial\Omega \leq \delta$, then there exists a constant $a$ such that
\begin{equation*}\label{e.l4.0}
  \|u-ax_n\|_{L^{\infty}(\Omega_{\eta})}\leq \eta^{1+\alpha}
\end{equation*}
and
\begin{equation*}\label{e.14.2}
|a|\leq C_1,
\end{equation*}
where $\eta$ depends only on $n,\lambda,\Lambda$ and $\alpha$.
\end{lemma}

\proof We prove the lemma by contradiction. Suppose that the lemma is false. Then there exist $0<\alpha<\alpha_1$ and sequences of $u_k,f_k,g_k,\Omega_k$ such that
\begin{equation*}
\left\{\begin{aligned}
&u_k\in S(\lambda,\Lambda,f_k)&& ~~\mbox{in}~~\Omega_k\cap B_1;\\
&u_k=g_k&& ~~\mbox{on}~~\partial \Omega_k\cap B_1
\end{aligned}\right.
\end{equation*}
with $\|u_k\|_{L^{\infty}(\Omega_k\cap B_1)}\leq 1$,
$\|f_k\|_{L^{n}(\Omega_k\cap B_1)}\leq 1/k$, $\|g_k\|_{L^{\infty}(\partial \Omega_k\cap B_1)}\leq 1/k$
and $\underset{B_1}{\mathrm{osc}}~\partial\Omega\leq 1/k$, and
\begin{equation}\label{e.l4.1}
  \|u_k-ax_n\|_{L^{\infty}(\Omega_{k}\cap B_{\eta})}> \eta^{1+\alpha}, \forall~|a|\leq C_1,
\end{equation}
where $0<\eta<1$ is taken small such that
\begin{equation}\label{e.2}
C_1\eta^{\alpha_1-\alpha}<1/2.
\end{equation}

Note that $u_k$ are uniformly bounded. In addition, by \Cref{l-7}, $u_k$ are equicontinuous. More precisely, for any $\Omega'\subset\subset B_1^+\cup T_1$, $\varepsilon>0$, there exist $\delta>0$ and $k_0$ such that for any $k\geq k_0$ and $x,y\in \Omega'\cap \bar{\Omega}_k$ with $|x-y|<\delta$, $|u(x)-u(y)|\leq \varepsilon$. Hence, there exists a subsequence (denoted by $u_k$ again) such that $u_k$ converges uniformly to some function $u$ on compact subsets of $B_1^+\cup T_1$. By the closedness (\Cref{l-5}), $u$ satisfies
\begin{equation*}
\left\{\begin{aligned}
&u\in S(\lambda,\Lambda,0)&& ~~\mbox{in}~~B_{1}^+;\\
&u=0&& ~~\mbox{on}~~T_{1}.
\end{aligned}\right.
\end{equation*}

By \Cref{l-1}, there exists $\bar{a}$ such that
\begin{equation*}
  |u(x)-\bar{a}x_n|\leq C_1 |x|^{1+\alpha_1}, ~~\forall ~x\in B_{1/2}^+
\end{equation*}
and
\begin{equation*}
  |\bar{a}|\leq C_1.
\end{equation*}
Hence, by noting\cref{e.2}, we have
\begin{equation}\label{e.3}
  \|u-\bar{a}x_n\|_{L^{\infty}(B_{\eta}^+)}\leq \eta^{1+\alpha}/2.
\end{equation}

By \Cref{l-7}, for $\delta$ small and $k$ large, we have
\begin{equation*}
    \|u_k-\bar{a}x_n\|_{L^{\infty}(\Omega_k\cap B_{\eta}\cap \left\{x_n\leq \delta\right\})}< \eta^{1+\alpha}.
\end{equation*}
Hence, from\cref{e.l4.1},
\begin{equation*}
  \|u_k-\bar{a}x_n\|_{L^{\infty}(\Omega_k\cap B_{\eta}\cap \left\{x_n>\delta\right\})}> \eta^{1+\alpha}.
\end{equation*}
Let $k\rightarrow \infty$, we have
\begin{equation*}
    \|u-\bar{a}x_n\|_{L^{\infty}(B_{\eta}\cap \left\{x_n>\delta\right\})}> \eta^{1+\alpha},
\end{equation*}
which contradicts with\cref{e.3}.  ~\qed~\\
\begin{remark}\label{r-7}
As pointed out in \cite[Chapte 1.3]{Wang_Regularity}, the benefits of the method of compactness are that it doesn't need the solvability of some equation, and the difference between the solution and the auxiliary function doesn't need to satisfy some equation.
\end{remark}

Now, we can prove the boundary $C^{1,\alpha}$ regularity.

\noindent\textbf{Proof of \Cref{t-1}.} We make some normalization first. Let $K_{\Omega}=[\partial \Omega]_{C^{1,\alpha}(0)}$. Then
\begin{equation}\label{e.1-1}
|x_n|\leq K_{\Omega}|x'|^{1+\alpha}, ~\forall ~x\in (\partial \Omega)_1.
\end{equation}

Next, we assume that $g(0)=0$ and $Dg(0)=0$. Otherwise, we may consider $v(x)=u(x)-g(0)-Dg(0)\cdot x$. Then the regularity of $u$ follows easily from that of $v$. Let $K_{g}=[g]_{C^{1,\alpha}(0)}$. Then
\begin{equation}\label{e.n3}
 |g(x)| \leq K_{g} |x|^{1+\alpha},~\forall~x\in (\partial \Omega)_1.
\end{equation}

Let $\delta$ be as in \Cref{l-4}. We assume that $\|u\|_{L^{\infty}(\Omega_1)}\leq 1$, $K_f\leq \delta$, $K_g\leq \delta/2$ and $K_{\Omega}\leq \delta/C_0$ where $C_0$ is a constant (depending only on $n,\lambda,\Lambda$ and $\alpha$) to be specified later. Otherwise, we may consider
\begin{equation*}
  v(y)=\frac{u(x)}{\|u\|_{L^{\infty}(\Omega_1)}+\delta^{-1}\left(K_f+2K_g\right)},
\end{equation*}
where $y=x/R$. By choosing $R$ small enough (depending only on $n,\lambda,\Lambda$ and $K_{\Omega}$), the above assumptions can be guaranteed. Without loss of generality, we assume that $R=1$.

To prove that $u$ is $C^{1,\alpha}$ at $0$, we only need to prove the following. There exists a sequence $a_k$ ($k\geq -1$) such that for all $k\geq 0$

\begin{equation}\label{e1.16}
\|u-a_kx_n\|_{L^{\infty }(\Omega _{\eta^{k}})}\leq \eta ^{k(1+\alpha )}
\end{equation}
and
\begin{equation}\label{e1.17}
|a_k-a_{k-1}|\leq C_1\eta ^{k\alpha},
\end{equation}
where $C_1$ is the universal constant as in \Cref{l-1} and $\eta$, depending only on $n,\lambda,\Lambda$ and $\alpha$, is as in \Cref{l-4} .

We prove the above by induction. For $k=0$, by setting $a_0=a_{-1}=0$, the conclusion holds clearly. Suppose that the conclusion holds for $k=k_0$. We need to prove that the conclusion holds for $k=k_0+1$.

Let $r=\eta ^{k_{0}}$, $y=x/r$ and
\begin{equation}\label{e-v1}
  v(y)=\frac{u(x)-a_{k_0}x_n}{r^{1+\alpha}}.
\end{equation}
Then $v$ satisfies
\begin{equation*}
\left\{\begin{aligned}
&v\in S(\lambda,\Lambda,\tilde{f})&& ~~\mbox{in}~~\tilde{\Omega}\cap B_1;\\
&v=\tilde{g}&& ~~\mbox{on}~~\partial \tilde{\Omega}\cap B_1,
\end{aligned}\right.
\end{equation*}
where
\begin{equation*}
  \tilde{f}(y)=\frac{f(x)}{r^{\alpha-1}},~\tilde{g}(y)=\frac{g(x)-a_{k_0}x_n}{r^{1+\alpha}}
  ~~\mbox{and}~~  \tilde{\Omega}=\frac{\Omega}{r}.
\end{equation*}

By\cref{e1.17}, there exists a constant $C_0$ depending only on $n,\lambda,\Lambda$ and $\alpha$ such that $|a_{k}|\leq C_0/2$ ($\forall~0\leq k\leq k_0$). Then it is easy to verify that
\begin{equation*}
  \|v\|_{L^{\infty}(\tilde{\Omega}\cap B_1)}\leq 1, ~(\mathrm{by}\cref{e1.16} ~\mathrm{and}\cref{e-v1})
\end{equation*}
\begin{equation*}
  \|\tilde{f}\|_{L^{n}(\tilde{\Omega}\cap B_1)}=\frac{\|f\|_{L^{n}(\Omega\cap B_r)}}{r^{\alpha}}\leq K_f\leq \delta, ~(\mathrm{by}\cref{e-dini-f})
\end{equation*}
\begin{equation}\label{e.g}
  \|\tilde{g}\|_{L^{\infty}(\partial \tilde{\Omega}\cap B_1)}\leq \frac{1}{r^{1+\alpha}}\left(K_gr^{1+\alpha}+\frac{C_0K_{\Omega}r^{1+\alpha}}{2}\right)\leq \delta  ~(\mathrm{by}\cref{e.1-1} ~\mathrm{and}\cref{e.n3})
\end{equation}
and
\begin{equation*}
\underset{B_1}{\mathrm{osc}}~\partial\tilde{\Omega}=
\frac{1}{r}\underset{B_r}{\mathrm{osc}}~\partial\Omega \leq K_{\Omega} r^{\alpha} \leq \delta.
\end{equation*}

By \Cref{l-4}, there exists a constant $\bar{a}$ such that
\begin{equation*}
\begin{aligned}
    \|v-\bar{a}y_n\|_{L^{\infty }(\tilde{\Omega} _{\eta})}&\leq \eta ^{1+\alpha}
\end{aligned}
\end{equation*}
and
\begin{equation*}
|\bar{a}|\leq C_1.
\end{equation*}

Let $a_{k_0+1}=a_{k_0}+r^{\alpha}\bar{a}$. Then\cref{e1.17} holds for $k_0+1$. Recalling\cref{e-v1}, we have
\begin{equation*}
  \begin{aligned}
&\|u-a_{k_0+1}x_n\|_{L^{\infty}(\Omega_{\eta^{k_0+1}})}\\
&= \|u-a_{k_0}x_n-r^{\alpha}\bar{a}x_n\|_{L^{\infty}(\Omega_{\eta r})}\\
&= \|r^{1+\alpha}v-r^{1+\alpha}\bar{a}y_n\|_{L^{\infty}(\tilde{\Omega}_{\eta})}\\
&\leq r^{1+\alpha}\eta^{1+\alpha}=\eta^{(k_0+1)(1+\alpha)}.
  \end{aligned}
\end{equation*}
Hence,\cref{e1.16} holds for $k=k_0+1$. By induction, the proof is completed.\qed~\\

\begin{remark}\label{r-5}
From the above proof, it shows clearly that the reason for the requirement of $\partial \Omega \in C^{1,\alpha}(0)$ is to estimate $x_n$ on $\partial \Omega$ (see\cref{e.g}). This observation is originated from \cite{MR3780142} and is key to the $C^{2,\alpha}$ regularity below.
\end{remark}

\section{$C^{2,\alpha}$ regularity}\label{S4}
In the following, we prove the boundary $C^{2,\alpha}$ regularity. From the proof for the $C^{1,\alpha}$ regularity, it can be inferred that if
\begin{equation*}
  \underset{B_r}{\mathrm{osc}}~\partial\Omega\leq Cr^{2+\alpha}, ~\forall ~0<r<1,
\end{equation*}
the $C^{2,\alpha}$ regularity follows almost exactly as the $C^{1,\alpha}$ regularity. However, the above can't be guaranteed by choosing a proper coordinate system, which is different from the $C^{1,\alpha}$ regularity. As pointed above, the requirement for $\partial \Omega$ is to estimate $x_n$ on $\partial \Omega$. If the term $x_n$ vanishes, the requirement for $\partial \Omega$ may be relaxed. This is the key idea for the $C^{2,\alpha}$ regularity.

The following lemma is similar to \Cref{l-4}, but without the term $x_n$ in the estimate.
\begin{lemma}\label{l-6}
Let $\alpha_2$ and $C_2$ be as in \Cref{l-2}. For any $0<\alpha<\alpha_2$, there exists $\delta>0$ such that if $u$ satisfies
\begin{equation*}
\left\{\begin{aligned}
&F(D^2u)=f&& ~~\mbox{in}~~\Omega_1;\\
&u=g&& ~~\mbox{on}~~(\partial \Omega)_1,
\end{aligned}\right.
\end{equation*}
with $\|u\|_{L^{\infty}(\Omega_1)}\leq 1$, $Du(0)=0$, $\|f\|_{L^{\infty}(\Omega_1)}\leq \delta$, $[g]_{C^{1,\alpha}(0)}\leq 1$, $\|g\|_{L^{\infty}((\partial \Omega)_1)}\leq \delta$ and $\|\partial \Omega\|_{C^{1,\alpha}(0)} \leq \delta$, then there exist constants $b_{in}$ such that
\begin{equation*}\label{e.l6.0}
  \|u-b_{in}x_ix_n\|_{L^{\infty}(\Omega_{\eta})}\leq \eta^{2+\alpha},
\end{equation*}
\begin{equation*}\label{e.l6.5}
F(b_{in})=0
\end{equation*}
and
\begin{equation*}\label{e.16.2}
|b_{in}|\leq C_2+1,
\end{equation*}
where $\eta$ depends only on $n,\lambda,\Lambda$ and $\alpha$.
\end{lemma}

\proof As before, we prove the lemma by contradiction. Suppose that the lemma is false. Then there exist $0<\alpha<\alpha_2$ and sequences of $F_k,u_k,f_k,g_k,\Omega_k$ such that
\begin{equation*}
\left\{\begin{aligned}
&F_k(D^2u_k)=f_k&& ~~\mbox{in}~~\Omega_k\cap B_1;\\
&u_k=g_k&& ~~\mbox{on}~~\partial \Omega_k\cap B_1,
\end{aligned}\right.
\end{equation*}
with $\|u_k\|_{L^{\infty}(\Omega_k\cap B_1)}\leq 1$, $[g_k]_{C^{1,\alpha}(0)}\leq 1$, $\|g_k\|_{L^{\infty}(\partial \Omega_k\cap B_1)}\leq 1/k$, $\|f_k\|_{L^{\infty}(\Omega_k\cap B_1)}\leq 1/k$, $\|\partial \Omega_k\|_{C^{1,\alpha}(0)} \leq 1/k$, $Du_k(0)=0$ and
\begin{equation}\label{e.l6.1}
  \|u_k-b_{in}x_ix_n\|_{L^{\infty}(\Omega_{k}\cap B_{\eta})}> \eta^{2+\alpha}, \forall~|b_{in}|\leq C_2+1 ~\mathrm{with}~F(b_{in})=0,
\end{equation}
where $0<\eta<1$ is taken small such that
\begin{equation}\label{e.l6.2}
C_2\eta^{\alpha_2-\alpha}<1/2.
\end{equation}

Since $F_k(0)=0$ and $F_k$ are Lipschitz continuous with a uniform Lipschitz constant depending only on $n,\lambda$ and $\Lambda$, there exists $F$ such that $F_k\rightarrow F$ on compact subsets of $S^{n}$. On the other hand, as before, $u_k$ are uniformly bounded and equicontinuous. Hence, by \Cref{l-5}, we can assume that $u_k$ converges uniformly to some function $u$ on compact subsets of $ B_1^+\cup T_1$ and $u$ satisfies
\begin{equation*}
\left\{\begin{aligned}
&F(D^2u)=0&& ~~\mbox{in}~~B_{1}^+;\\
&u=0&& ~~\mbox{on}~~T_{1}.
\end{aligned}\right.
\end{equation*}

By the $C^{1,\alpha}$ estimate for $u_k$ (see \Cref{t-1}) and noting $Du_k(0)=0$, we have
\begin{equation*}
\|u_k\|_{L^{\infty }(\Omega_k\cap B_r)}\leq Cr^{1+\bar{\alpha}} ~~~~\forall~0<r<1,
\end{equation*}
where $\bar{\alpha}<\min(\alpha,\alpha_1)$ and $C$ is universal. Since $u_k$ converges to $u$ uniformly,
\begin{equation*}
\|u\|_{L^{\infty }(B_r^+)}\leq Cr^{1+\bar{\alpha}} ~~~~\forall~0<r<1,
\end{equation*}
Hence, $Du(0)=0$.

By \Cref{l-2}, there exist $\bar{b}_{in}$ such that
\begin{equation*}
  |u(x)-\bar{b}_{in}x_ix_n|\leq C_2 |x|^{2+\alpha_2}, ~~\forall ~x\in B_{1/2}^+,
\end{equation*}
\begin{equation*}
  F(\bar{b}_{in})=0
\end{equation*}
and
\begin{equation*}
  |\bar{b}_{in}|\leq C_2.
\end{equation*}

Since $F_k(\bar{b}_{in})\rightarrow F(\bar{b}_{in})=0$. For $k$ large, there exists $t_k$ with $|t_k|\leq \eta^{2+\alpha}/4$ and $t_k\rightarrow 0$ such that
\begin{equation*}
  F_k(\bar{b}_{in}+t_k\delta_{nn})=0,
\end{equation*}
where $\delta_{nn}$ denotes the matrix $a_{ij}$ whose elements are all $0$ except $a_{nn}=1$ (similarly hereinafter).

By noting\cref{e.l6.2}, we have
\begin{equation}\label{e.l6.3}
  \|u-\bar{b}_{in}x_ix_n\|_{L^{\infty}(B_{\eta}^+)}\leq \eta^{2+\alpha}/2.
\end{equation}

By \Cref{l-7}, for $\delta$ small and $k$ large, we have
\begin{equation*}
    \|u_k-\bar{b}_{in}x_ix_n-t_kx_n^2\|_{L^{\infty}(\Omega_k\cap B_{\eta}\cap \left\{x_n\leq \delta\right\})}< \eta^{2+\alpha}.
\end{equation*}
Hence, from \cref{e.l6.1},
\begin{equation*}
  \|u_k-\bar{b}_{in}x_ix_n-t_kx_n^2\|_{L^{\infty}(\Omega_k\cap B_{\eta}\cap \left\{x_n>\delta\right\})}> \eta^{2+\alpha}.
\end{equation*}
Let $k\rightarrow \infty$, we have
\begin{equation*}
    \|u-\bar{b}_{in}x_ix_n\|_{L^{\infty}(B_{\eta}\cap \left\{x_n>\delta\right\})}> \eta^{2+\alpha},
\end{equation*}
which contradicts with\cref{e.l6.3}.  ~\qed~\\

The following is the essential result for the $C^{2,\alpha}$ regularity. The key is that if $Du(0)=0$, the $C^{2,\alpha}$ regularity holds even if $\partial \Omega\in C^{1,\alpha}(0)$.
\begin{theorem}\label{t-4}
Let $0<\alpha <\alpha _2$ and $\partial \Omega$ be $C^{1,\alpha}$ at $0$. Assume that $u$ satisfies
\begin{equation*}
\left\{\begin{aligned}
&F(D^2u)=f&& ~~\mbox{in}~~\Omega_1;\\
&u=g&& ~~\mbox{on}~~(\partial \Omega)_1,
\end{aligned}\right.
\end{equation*}
with $Du(0)=0$. Suppose that
\begin{equation}\label{e.4-3}
 |f(x)|\leq K_f|x|^{\alpha}, ~~\forall ~x\in (\partial \Omega)_1.
\end{equation}
and
\begin{equation}\label{e.4-4}
|g(x)|\leq K_g|x|^{2+\alpha}, ~~\forall ~x\in (\partial \Omega)_1.
\end{equation}

Then $u\in C^{2,\alpha}(0)$, i.e., there exists a quadratic polynomial $P$ such that
\begin{equation}\label{e.n1}
  |u(x)-P(x)|\leq C |x|^{2+\alpha}\left(\|u\|_{L^{\infty }(\Omega_1)}+K_f+K_g\right), ~~\forall ~x\in \Omega_{r_1},
\end{equation}
and
\begin{equation}\label{e.nt2-2}
|Du(0)|+\|D^2u(0)\|\leq C\left(\|u\|_{L^{\infty }(\Omega_1)}+K_f+K_g\right),
\end{equation}
where $C$ depends only on $n, \lambda, \Lambda$ and $\alpha$, and $r_1$ depends also on $[\partial \Omega]_{C^{1,\alpha}(0)}$.
\end{theorem}

\proof As before, we make some normalization first. Let $K_{\Omega}=[\partial \Omega]_{C^{1,\alpha}(0)}$. Then
\begin{equation}\label{e.n1-1}
|x_n|\leq K_{\Omega}|x'|^{1+\alpha}, ~\forall ~x\in (\partial \Omega)_1.
\end{equation}

Let $\delta$ be as in \Cref{l-6}. As before, we assume that $\|u\|_{L^{\infty}(\Omega_1)}\leq 1$, $K_f\leq \delta$, $K_g\leq \delta/2$ and $K_{\Omega}\leq \delta/C_0$ where $C_0$ is a constant (depending only on $n,\lambda,\Lambda$ and $\alpha$) to be specified later.

To prove that $u$ is $C^{2,\alpha}$ at $0$, we only need to prove the following. There exist sequences $ (b_k)_{in}$ ($k\geq -1$) such that for all $k\geq 0$,

\begin{equation}\label{e.t4.5}
\|u-(b_k)_{in}x_ix_n\|_{L^{\infty }(\Omega _{\eta^{k}})}\leq \eta ^{k(2+\alpha )},
\end{equation}
\begin{equation}\label{t.t4.6}
F((b_k)_{in})=0
\end{equation}
and
\begin{equation}\label{e.t4.7}
|(b_k)_{in}-(b_{k-1})_{in}|\leq (C_2+1)\eta ^{k\alpha},
\end{equation}
where $C_2$ is the universal constant as in \Cref{l-2} and $\eta$, depending only on $n,\lambda,\Lambda$ and $\alpha$, is as in \Cref{l-6}.

We prove the above by induction. For $k=0$, by setting $(b_0)_{in}=(b_{-1})_{in}=0$, the conclusion holds clearly. Suppose that the conclusion holds for $k=k_0$. We need to prove that the conclusion holds for $k=k_0+1$.

Let $r=\eta ^{k_{0}}$, $y=x/r$ and
\begin{equation}\label{e.t4.8}
  v(y)=\frac{u(x)-(b_{k_0})_{in}x_ix_n}{r^{2+\alpha}}.
\end{equation}
Then $v$ satisfies
\begin{equation*}
\left\{\begin{aligned}
&\tilde{F}(D^2v)=\tilde{f}&& ~~\mbox{in}~~\tilde{\Omega}\cap B_1;\\
&v=\tilde{g}&& ~~\mbox{on}~~\partial \tilde{\Omega}\cap B_1,
\end{aligned}\right.
\end{equation*}
where for $M\in S^{n\times n}$,
\begin{equation*}
\tilde{F}(M)=\frac{F(r^{\alpha}M+(b_{k_0})_{in})}{r^{\alpha}},
 \tilde{f}(y)=\frac{f(x)}{r^{\alpha}},
 \tilde{g}(y)=\frac{g(x)-(b_{k_0})_{in}x_ix_n}{r^{2+\alpha}}
  ~~\mbox{and}~~  \tilde{\Omega}=\frac{\Omega}{r}.
\end{equation*}

Then $\tilde{F}$ is uniformly elliptic with ellipticity constants $\lambda$ and $\Lambda$ and $\tilde{F}(0)=0$. By\cref{e.t4.7}, there exists a constant $C_0$ depending only on $n,\lambda,\Lambda$ and $\alpha$ such that $|(b_k)_{in}|\leq C_0/2$ ($\forall~0\leq k\leq k_0$). Then it is easy to verify that
\begin{equation*}
  \|v\|_{L^{\infty}(\tilde{\Omega}\cap B_1)}\leq 1, ~(\mathrm{by}\cref{e.t4.5} ~\mbox{and}\cref{e.t4.8})
\end{equation*}
\begin{equation*}
  \|\tilde{f}\|_{L^{\infty}(\tilde{\Omega}\cap B_1)}=\frac{\|f\|_{L^{\infty}(\Omega\cap B_r)}}{r^{\alpha}}\leq K_f\leq \delta ~(\mathrm{by}\cref{e.4-3})
\end{equation*}
and
\begin{equation*}
\|\partial \tilde{\Omega}\cap B_1\|_{C^{1,\alpha}(0)} \leq K_{\Omega} r^{\alpha} \leq \delta.
\end{equation*}
In addition, by\cref{e.4-4} and\cref{e.n1-1}, we have
\begin{equation*}
  \begin{aligned}
      \|\tilde{g}\|_{L^{\infty}(\partial \tilde{\Omega}\cap B_t)} &\leq \frac{1}{r^{2+\alpha}}\left(K_gt^{2+\alpha}r^{2+\alpha}+
\frac{C_0K_{\Omega}t^2r^{2+\alpha}}{2}\right)\leq \delta t^2.
  \end{aligned}
\end{equation*}
Hence,
\begin{equation*}
  [\tilde{g}]_{C^{1,\alpha}(0)}\leq \delta \leq 1 ~\mbox{and}~\|\tilde{g}\|_{L^{\infty}(\partial \tilde{\Omega}\cap B_1)}\leq \delta.
\end{equation*}

By \Cref{l-6}, there exists constants $\bar{b}_{in}$ such that
\begin{equation*}
\begin{aligned}
    \|v-\bar{b}_{in}y_iy_n\|_{L^{\infty }(\tilde{\Omega} _{\eta})}&\leq \eta ^{2+\alpha},
\end{aligned}
\end{equation*}
\begin{equation*}
  \tilde{F}(\bar{b}_{in})=0
\end{equation*}
and
\begin{equation*}\label{e.19}
|\bar{b}_{in}|\leq C_2+1.
\end{equation*}

Let $(b_{k_0+1})_{in}=(b_{k_0})_{in}+r^{\alpha}\bar{b}_{in}$. Then\cref{t.t4.6} and\cref{e.t4.7} hold for $k_0+1$. Recalling\cref{e.t4.8}, we have
\begin{equation*}
  \begin{aligned}
&\|u-(b_{k_0+1})_{in}x_ix_n\|_{L^{\infty}(\Omega_{\eta^{k_0+1}})}\\
&= \|u-(b_{k_0})_{in}x_ix_n-r^{\alpha}\bar{b}_{in}x_ix_n\|_{L^{\infty}(\Omega_{\eta r})}\\
&= \|r^{2+\alpha}v-r^{2+\alpha}\bar{b}_{in}y_iy_n\|_{L^{\infty}(\tilde{\Omega}_{\eta})}\\
&\leq r^{2+\alpha}\eta^{2+\alpha}=\eta^{(k_0+1)(2+\alpha)}.
  \end{aligned}
\end{equation*}
Hence,\cref{e.t4.5} holds for $k=k_0+1$. By induction, the proof is completed.\qed~\\

\noindent\textbf{Proof of \Cref{t-2}.} In fact, \Cref{t-4} has contained the essential ingredients for the $C^{2,\alpha}$ regularity. The following proof is just the normalization in some sense.

Assume that $\Omega$ satisfies\cref{e-re} and\cref{e-re2} with $P(x')=x'^{T}Ax'$ for some $A\in S^{n\times n}$. By scaling, we can assume that $\|\partial \Omega\|_{C^{2,\alpha}(0)}\leq 1$.

Let $F_1(M)=F(M)-f(0)$ for $M\in S^{n\times n}$. (In the following proof, $M$ always denotes a symmetric matrix.) Then $F_1$ is uniformly elliptic with the same ellipticity constants and $u$ satisfies
\begin{equation*}
\left\{\begin{aligned}
&F_1(D^2u)=f_1&& ~~\mbox{in}~~\Omega_1;\\
&u=g&& ~~\mbox{on}~~(\partial \Omega)_1,
\end{aligned}\right.
\end{equation*}
where $f_1(x)=f(x)-f(0)$.

Next, let $u_1(x)=u(x)-g(0)-Dg(0)\cdot x-x^TD^2g(0)x$ and $F_2(M)=F_1(M+D^2g(0))$. Then $F_2$ is uniformly elliptic with the same ellipticity constants and $u_1$ satisfies
\begin{equation*}
\left\{\begin{aligned}
&F_2(D^2u_1)=f_1&& ~~\mbox{in}~~\Omega_1;\\
&u_1=g_1&& ~~\mbox{on}~~(\partial \Omega)_1,
\end{aligned}\right.
\end{equation*}
where $g_1(x)=g(x)-g(0)-Dg(0)\cdot x-x^TD^2g(0)x$. Hence,
\begin{equation*}
  |f_1(x)|\leq [f]_{C^{\alpha}(0)}|x|^{\alpha}, ~~\forall ~x\in \Omega_1,
\end{equation*}
\begin{equation*}
  |g_1(x)|\leq [g]_{C^{2,\alpha}(0)}|x|^{2+\alpha}, ~~\forall ~x\in (\partial \Omega)_1
\end{equation*}
and
\begin{equation*}
  |F_2(0)|=|F_1(D^2g(0))|=|F(D^2g(0))-f(0)|\leq C\left(\|D^2g(0)\|+|f(0)|\right),
\end{equation*}
where $C$ is universal.

Note that (see \cite[Proposition 2.13]{MR1351007}),
\begin{equation*}
  u_1\in S(\lambda/n,\Lambda,f_1-F_2(0)).
\end{equation*}
Then by \Cref{t-1}, $u_1\in C^{1,\bar{\alpha}}(0)$ for $\bar{\alpha}=\min(\alpha_1,\alpha_2)/2$, $Du_1(0)=(0,...,0,(u_1)_n(0))$ and
\begin{equation}\label{e.4-5}
  \begin{aligned}
|(u_1)_n(0)| &\leq C \left(\|u_1\|_{L^{\infty }(\Omega\cap B_1)}+[f]_{C^{\alpha}(0)}+|F_2(0)|+[g]_{C^{2,\alpha}(0)}\right)\\
    &\leq C\left(\|u\|_{L^{\infty }(\Omega\cap B_1)}+\|f\|_{C^{\alpha}(0)}+\|g\|_{C^{2,\alpha}(0)}+|F_2(0)|\right)\\
    &\leq C\left(\|u\|_{L^{\infty }(\Omega\cap B_1)}+\|f\|_{C^{\alpha}(0)}+\|g\|_{C^{2,\alpha}(0)}\right),
  \end{aligned}
\end{equation}
where $C$ is universal.

Let $u_2(x)=u_1(x)-(u_1)_n(0)\left(x_n-x'^{T}Ax'\right)$ and $F_3(M)=F_2(M-(u_1)_n(0)A)$. Then $F_3$ is uniformly elliptic with the same ellipticity constants and $u_2$ satisfies
\begin{equation*}
\left\{\begin{aligned}
&F_3(D^2u_2)=f_1&& ~~\mbox{in}~~\Omega_1;\\
&u_2=g_2&& ~~\mbox{on}~~(\partial \Omega)_1,
\end{aligned}\right.
\end{equation*}
where $g_2=g_1-(u_1)_n(0)\left(x_n-x'^{T}Ax'\right)$.

Next, let $u_3(x)=u_2(x)+tx_n^2$ and $F_4(M)=F_3(M-2t\delta_{nn})$. Then $F_4(0)=0$ for some $t\in R$ and (note that $\|A\|\leq \|\partial \Omega\|_{C^{2,\alpha}(0)}\leq 1$)
\begin{equation}\label{e.4-6}
  \begin{aligned}
|t|&\leq |F_3(0)|/\lambda\leq C |F(D^2g(0)-(u_1)_n(0)A)-f(0)|\\
    &\leq C\left(\|D^2g(0)\|+|(u_1)_n(0)|\|A\|+|f(0)|\right)\\
    &\leq C \left(\|u\|_{L^{\infty }(\Omega\cap B_1)}+\|f\|_{C^{\alpha}(0)}+\|g\|_{C^{2,\alpha}(0)}\right),
  \end{aligned}
\end{equation}
where $C$ is universal. Moreover, $u_3$ satisfies
\begin{equation*}
\left\{\begin{aligned}
&F_4(D^2u_3)=f_1&& ~~\mbox{in}~~\Omega_1;\\
&u_3=g_3&& ~~\mbox{on}~~(\partial \Omega)_1,
\end{aligned}\right.
\end{equation*}
where $g_3=g_2+tx_n^2$.

Then it is easy to verify that $F_4(0)=0$, $Du_3(0)=0$ and
\begin{equation*}
  \begin{aligned}
    |g_3(x)|& \leq \left([g]_{C^{2,\alpha}(0)}+|(u_1)_n(0)|[\partial \Omega]_{C^{2,\alpha}(0)}+|t|\|\partial \Omega\|^2_{C^{2,\alpha}(0)}\right)|x|^{2+\alpha}\\
    &\leq C\left(\|u\|_{L^{\infty }(\Omega\cap B_1)}+\|f\|_{C^{\alpha}(0)}+\|g\|_{C^{2,\alpha}(0)}\right)|x|^{2+\alpha}, ~~\forall ~x\in (\partial \Omega)_1,
  \end{aligned}
\end{equation*}
where $C$ is universal.

By \Cref{t-4}, $u_3$ and hence $u$ is $C^{2,\alpha}$ at $0$, and the estimates\cref{en2} and\cref{e.t2-2} hold. \qed~\\

\section*{Acknowledgments}
The authors would like to thank Professor Dongsheng Li for useful discussions.

\section*{References}
\bibliographystyle{model4-names}

\end{document}